\documentclass[11pt]{amsart}
\usepackage[margin=2.6cm]{geometry}
\usepackage{amssymb,amsmath}
\newtheorem{theorem}{Theorem}[section]
\newtheorem{definition}[theorem]{Definition}
\def \r{{\bf H}_{\R}}
\def\z{{\bf z}}
\def\w{{\bf w}}
\def\P{\mathbb P}
\def\R{\mathbb R}
\def\V{\mathrm V}
\def\C{\mathbb C}
\def \Ad{\mathrm{Ad}}
\def \c{{\bf H}_{\C}}

\def\PU{\mathrm{PU}}

\newtheorem{corollary}[theorem]{Corollary}
\newtheorem{proposition}[theorem]{Proposition}
\newtheorem{lemma}[theorem]{Lemma}

\def\r{\mathcal R}
\def\R{\mathbb R}

\def\C{\mathbb C}

\def\H{\mathbb H}

\def\z{{\bf z}}
\def\w{{\bf w}}
\def\a{{\bf a}}

\def\V{\mathrm V}

\def \h{{\bf H}_{\H}}
\def\R{\mathbb R}

\newcommand{\SO}{\mathrm{SO}}

\newcommand{\Sp}{\mathrm{Sp}}
\newcommand{\PSp}{\mathrm{PSp}}

\def\P{\mathbb P}

\def \a {{\bf a}}
\def \r  {{\bf r}}
\def \x {{\bf x}}

\newcommand{\SU}{\mathrm{SU}}
\newtheorem{prop}[theorem]{Proposition}
\theoremstyle{definition}

\theoremstyle{remark}
\newtheorem{remark}{Remark}
\numberwithin{equation}{section}

\numberwithin{equation}{section}

\begin{document}
\title[Strongly Doubly Reversible Pairs in $\PSp(n,1)$]{Strongly Doubly Reversible  Pairs in Quaternionic Unitary Group of Signature $(n,1)$}
\author[Krishnendu Gongopadhyay  \and Sagar B. Kalane]{Krishnendu Gongopadhyay  \and Sagar B. Kalane}

\address{Indian Institute of Science Education and Research (IISER) Mohali,
	Knowledge City,  Sector 81, S.A.S. Nagar 140306, Punjab, India}
\email{krishnendu@iisermohali.ac.in}
\address{Institute of Mathematical Sciences, IV Cross Road, CIT Campus Taramani, Chennai 600 113 Tamil Nadu, India.}
\email{sagarbk@imsc.res.in, sagark327@gmail.com}
\date{\today}
	\begin{abstract}
Let $\PSp(n,1)$ denote the isometry group of quaternionic hyperbolic $n$--space $\h^n$.
A pair $(g_1,g_2)\in\PSp(n,1)^2$ is \emph{strongly doubly reversible} if
$(g_1,g_2)$ and $(g_1^{-1},g_2^{-1})$ are simultaneously conjugate by an involution.
Equivalently, there exist involutions $i_1,i_2,i_3\in\PSp(n,1)$ such that
$g_1=i_1i_2$ and $g_2=i_1i_3$.

We prove that the set of strongly doubly reversible pairs has Haar measure zero in
$\PSp(n,1)\times\PSp(n,1)$. The same conclusion holds for
$\PSp(n)\times\PSp(n)$ and $\SU(n,1)\times\SU(n,1)$ for $n\ge2$, and for
$\SO_0(n,1)\times\SO_0(n,1)$ and $\SO(n+1)\times\SO(n+1)$ for $n\ge4$.
In the compact rank-one case, every pair in $\PSp(1)$ is strongly doubly
reversible, which gives a short proof of the theorem of Basmajian and Maskit
that every pair in ${\rm SO}(4)$ is strongly doubly reversible.

For hyperbolic pairs, we prove that double reversibility and strong double
reversibility are equivalent in $\PSp(1,1)$. In higher dimension, we prove the
same implication when one member of the pair is regular, with pairwise distinct
non-real unit eigenvalue classes. Finally, in $\PSp(1,1)$, for a hyperbolic element in normal form, we give a complete explicit matrix criterion characterizing all elements that form a strongly doubly reversible pair with it.
\end{abstract} 	
    \maketitle

\section{Introduction}

An element in a group $G$ is said to be \emph{strongly reversible} (or 
\emph{strongly real}) if it can be written as a product of two involutions in 
$G$. This notion is closely related to that of \emph{reversible} (or 
\emph{real}) elements, which are conjugate to their inverses in $G$. Every 
strongly reversible element is necessarily reversible, but the converse is not 
true, in general. The classification and structure of such elements have been 
the subject of regular investigation in various branches of mathematics, for 
example, see \cite{BR}, \cite{Dj}, \cite{El}, \cite{gm2}, \cite{os}, \cite{Sa}, \cite{ST}, 
\cite{Wo}. 

Beyond their algebraic significance, strongly reversible elements play a 
central role in understanding symmetries in geometry. In particular, certain 
geometrically natural groups are built entirely from such elements. A 
classical and striking example arises in the setting of hyperbolic geometry. 
The group $ \mathrm{PSL}(2, \mathbb{C}) $ can be identified with the 
orientation-preserving isometries of the three dimensional hyperbolic space. 
If  $ A $ and $ B $ are elements of $ \mathrm{PSL}(2, \mathbb{C}) $ 
generating a non-elementary subgroup, then there exist involutions $ i_1, 
i_2, i_3 \in \mathrm{PSL}(2, \mathbb{C}) $ such that $ A = i_1 i_2 $ and $ 
B = i_1 i_3 $, see ~\cite{gold}, \cite{m}. In the real hyperbolic case, a 
similar statement holds for pairs in $ \mathrm{PSL}(2, \mathbb{R}) $ with the 
involutions given by  orientation-reversing reflections. 

These observations motivate a more general concept that extends beyond 
individual decompositions to relationships between two elements and can be 
formulated for any abstract group. 
\begin{definition}
Let $ G $ be a group. Consider the $G$ action on $G \times G$ by 
simultaneous conjugation:
$$g.(g_1, g_2)=(gg_1g^{-1}, gg_2g^{-1}). $$
For two elements $ g_1, g_2 \in G$,  the pair $(g_1, g_2)$ is said to be 
\emph{doubly  reversible} or \emph{doubly real} if $(g_1, g_2)$ and 
$(g_1^{-1}, g_2^{-1})$ belong to the same $G$-conjugation orbit. Furthermore, 
if we choose the conjugating element $g$ to be such that it is an involution, 
i.e. $g^2=1$, then we call $(g_1, g_2)$ to be \emph{strongly doubly 
reversible} or \emph{strongly doubly real}. 
\end{definition}

 Strong double reversibility encodes the existence of a common involutory symmetry for a pair of elements, a condition that naturally arises in the study of symmetric representations of free groups and subgroups generated by involutions. It strengthens classical reversibility by requiring such a symmetry to be shared by two elements simultaneously. This condition is not merely technical as it captures the extent to which pairs of transformations are governed by involutory symmetries. A basic problem, therefore, is to understand whether such symmetric representations are generic or exceptional within a given group.

The above interpretation captures the geometric compatibility between  
elements $ g_1$ and $ g_2 $.  For a strongly doubly reversible pair $(g_1, 
g_2)$, both elements are generated by pairs of involutions that share a common 
factor.  The classical result for $ \mathrm{PSL}(2,\mathbb{C}) $ may thus be 
interpreted as asserting that any two generators of a non-elementary subgroup 
are necessarily  strongly doubly reversible. 

This notion can be extended to any $k$-tuple of elements in $G$ to define $k$-
reversible (or $k$-real) and strongly $k$-reversible (or strongly $k$-real)
tuples. Our main focus is the case $k=2$, although the measure-theoretic
argument of Section~\ref{spn1s} also gives a statement for $k$-tuples.

Note that if $(g_1, g_2)$ is strongly doubly reversible, then there exist 
involutions $ i_1, i_2, i_3 \in G $ such that $ g_1 = i_1 i_2 $ and $ g_2 
= i_1 i_3 $. Conversely, if there are involutions $i_1, i_2, i_3 \in G$ such 
that $ g_1 = i_1 i_2 $ and $ g_2 = i_1 i_3 $, then $(g_1, g_2)$ is 
strongly doubly reversible. In particular, every strongly doubly reversible 
pair is doubly reversible.

The study of such  pairs is particularly interesting in groups where every 
element is a product of two involutions. In such settings, one may naturally 
ask which pairs of elements in such a group are  strongly doubly reversible. 
This question not only provides understanding about the group's internal 
structure, but also connects it to broader topics like discreteness and 
geometric finiteness. For instance, when $ g_1 $ and $ g_2 $ are  strongly 
doubly reversible, the subgroup $ \langle g_1, g_2 \rangle $ sits as a subgroup of index at most two, of the group $ \langle i_1, i_2, i_3 \rangle $ generated 
by involutions. This can potentially lead to better insight into groups 
generated by three involutions, e.g., triangle groups in hyperbolic geometry. 
Despite its relevance, the problem of classifying  strongly doubly reversible 
elements remains largely unexplored. Even within the context of finite groups, 
systematic efforts to understand doubly reversible pairs have begun only 
recently, e.g. \cite{amri}. The terminology `$k$-real' has been borrowed from 
\cite{amri}. 

In geometric contexts, Basmajian and Maskit in~\cite{bm}, posed the 
problem of generalizing the classical $ \mathrm{PSL}(2,\mathbb{C})$ result 
to higher-dimensional M\"obius groups and isometry groups of Riemannian space 
forms.  It may be noted that strongly doubly reversible pairs were termed 
\emph{linked pairs} in \cite{bm}. Basmajian and Maskit proved that for higher 
dimensions, especially $n \geq 5$, almost all pairs in these groups are not 
strongly doubly reversible.  Basmajian and Maskit also proved that every pair 
in the orthogonal group ${\rm SO}(4)$ is strongly doubly reversible.   In a 
subsequent work, Silverio~\cite{S} provided a geometric description of 
strongly doubly reversible pairs in real hyperbolic $4$-space. In complex 
hyperbolic geometry, the strongly doubly reversible pairs acquire additional 
structure. When every element in ${\mathrm{PU}}(n,1)$, isometry group of the 
$n$-dimensional complex hyperbolic space, is a product of two anti-holomorphic involutions, not every element of ${\PU}(n,1)$ is a product of two 
(holomorphic) involutions, cf. \cite{gp2}. In $ \mathrm{PU}(2,1) $ they are 
named as \emph{$ \mathbb{R}$-decomposable} or \emph{$ \mathbb{C}$-
decomposable}, depending on whether the generating involutions are anti-
holomorphic or holomorphic. Will~\cite{w} classified the loxodromic pairs, 
while Paupert and Will~\cite{pw} gave a complete classification of the $ 
\mathbb{R}$-decomposable pairs in ${\rm PU}(2, 1)$. The $ \mathbb{C}$-
decomposable pairs in $ \mathrm{PU}(2,1) $ have been described by Ren et 
al.~\cite{ren}. 

It should be noted that no analogous classification for strongly doubly
reversible pairs is known for $\mathrm{PU}(n,1)$ with $n\ge3$,
$\mathrm{SO}_0(n,1)$ with $n\ge5$, or for most other linear groups,
including compact ones. Here $\mathrm{SO}_0(n,1)$ denotes the identity component of the special
orthogonal group associated to a quadratic form of signature $(n,1)$, and its quotient by $\{\pm I\}$ is identified with the group of
orientation-preserving isometries of real hyperbolic $n$-space.  For $\mathrm{SO}_0(n,1)$, $n\ge5$, the best result in
this direction appears in \cite{bm}, where a measure-theoretic observation
(counterpart of Theorem~\ref{spn1} of this paper) was obtained.

 Let $ \h^n $ denote the $n$-dimensional quaternionic hyperbolic space, whose 
isometry group is $ \PSp(n,1) = \Sp(n,1)/\{ \pm I \} $. A result by Bhunia 
and Gongopadhyay~\cite{bg} shows that every element of $ \Sp(n,1) $ can be 
expressed as a product of two \emph{skew-involutions}. Recall that a skew-
involution is an element $ i \in \Sp(n,1) $ such that $ i^2 = -1 $. The 
skew-involutions project to involutions in $ \PSp(n,1) $. In contrast to $ 
\PSp(n,1) $, the group $ \Sp(n,1) $ itself has relatively few genuine 
involutions, and not all elements can be written as products of two such. 
Since every element of $\PSp(n,1)$ is strongly reversible, it is natural to 
ask about strongly doubly reversible pairs in $\PSp(n,1)$.

In this paper, we do not attempt a complete 
classification of such pairs, as this seems out of reach 
due to the non-commutative nature of quaternionic 
multiplication and the inherent complexity of higher-dimensional geometry. 
The absence of a well-behaved 
trace function or comprehensive conjugacy invariants in 
the quaternionic setting further complicates the 
analysis. These invariants are, however, central to the 
classification of strongly doubly reversible pairs in 
$\mathrm{PU}(2,1)$, the isometry group of the two-dimensional complex hyperbolic space, see \cite{pw}. Accordingly, many standard tools from complex  hyperbolic geometry do not extend even to $\PSp(2,1)$ for understanding the pairs. 

Nevertheless, we observe some interesting results for strongly doubly 
reversible pairs in $\PSp(n,1)$. Now we outline the main results. The first is 
the following generalization of \cite[Theorem 1.5]{bm}. This shows that strongly doubly reversible pairs are non-generic in $\PSp(n, 1)$. 
\begin{theorem}\label{spn1}
 The set of strongly doubly reversible pairs in $\PSp(n,1)$  has Haar measure 
 zero in $\PSp(n, 1)\times \PSp(n,1)$. 
\end{theorem}

The above result shows that strong double reversibility is genuinely non-generic in $\PSp(n,1)$.  In particular, a random pair of elements in 
$\PSp(n,1)$ almost surely does not admit any common involutory symmetry, despite the fact that every individual element is strongly reversible.

  We prove the above result using an approach 
  based on Lie-theoretic considerations. This is in contrast to the methods of \cite{bm} for ${\rm SO}(n,1)$. While the 
  approach in \cite{bm} relies on a detailed geometric 
  analysis of `common perpendicular subspaces', our 
  method is instead based on structural properties of 
  involutions and dimension-count arguments. This perspective is uniform and representation-independent, and may apply in other settings as well. For instance, using similar arguments, the result extends to strongly doubly reversible pairs in $\SU(n,1)$ and $\PSp(n)$ for $n\ge2$, and gives an independent proof of the measure-theoretic result of Basmajian--Maskit for $\SO(n,1)$. 

  \begin{corollary}
   Let $n \geq 2$. The set of strongly doubly reversible pairs in $\SU(n,1)$
   has Haar measure zero in $\SU(n, 1)\times \SU(n,1)$. 
\end{corollary}
\begin{corollary}
  Let $n \geq 2$. The set of strongly doubly reversible pairs in $\PSp(n)$ has 
  Haar measure zero in $\PSp(n)\times \PSp(n)$. 
\end{corollary}

\begin{corollary}\label{cor:SO}
Let $n\ge4$. The set of strongly doubly reversible pairs in $\SO_0(n,1)$
has Haar measure zero in $\SO_0(n,1)\times\SO_0(n,1)$. The same conclusion holds
for $\SO(n+1)$.
\end{corollary}

However, when $n=1$, we see that every element in $\PSp(1)$ is  strongly
doubly reversible. We have used elementary quaternionic analysis to see this for  $\PSp(1)$.  We also apply this result to offer a very short proof of \cite[Theorem 1.4]{bm}, which is the following. 
\begin{theorem}
    Every pair of elements in ${\rm SO}(4)$ is strongly doubly reversible. 
\end{theorem}

Following the terminology in \cite{Rat}, recall that an element $g$ in
$\Sp(n,1)$ is called \emph{hyperbolic} if it has exactly two fixed points on the
boundary of $\h^n$. A hyperbolic element is semisimple, and up to conjugacy it is
represented by a complex diagonal matrix with two null eigenvalue classes represented by
$re^{i\theta}$ and $r^{-1}e^{i\theta}$, $0<r<1$, and the remaining eigenvalue classes
represented by unit complex numbers. A hyperbolic element is called \emph{regular} if $\theta \neq 0, \pi/2, \pi$,  and it has mutually distinct non-real unit eigenvalue classes. 

The following theorem relates double reversibility and strong double reversibility when one element of the pair is hyperbolic.
\begin{theorem}\label{hyp}
Let $A \in \PSp(n,1)$ be hyperbolic and let $B \in \PSp(n,1)$ be arbitrary.
\begin{enumerate}
\item If $n=1$, then $(A,B)$ is doubly reversible if and only if it is strongly
      doubly reversible.
\item If $n\ge 2$ and $A$ is regular, then double reversibility implies strong
      double reversibility.
\end{enumerate}
\end{theorem}
The regularity hypothesis in the second assertion is exactly what is needed for the blockwise argument used in the proof. Here, the null rotation angle is non-real and different from \(\pi/2\), and the unit eigenvalue classes are non-real and pairwise distinct. When these conditions fail, reversers may have additional freedom and need not themselves be involutions, so the above argument no longer applies. We do not address here whether double reversibility implies strong double reversibility for arbitrary hyperbolic pairs in higher dimensions.

For $\PSp(1,1)$ we also obtain a complete explicit matrix criterion. After conjugating a hyperbolic element to normal form, we characterize all elements of $\PSp(1,1)$ that form a strongly doubly reversible pair with it; see Theorem~\ref{qd} and Corollary~\ref{qd-nondegenerate}.

\subsubsection*{Structure of the paper}
After discussing notation and preliminaries in Section~\ref{pre}, we prove in
Section~\ref{sp11} that every pair of elements in $\PSp(1)$ is strongly doubly
reversible and derive the corresponding ${\rm SO}(4)$ application. In
Section~\ref{spn1s}, we establish the Haar-measure-zero theorem. Section~\ref{sp-1}
examines double versus strong double reversibility for hyperbolic pairs, first in
$\PSp(1,1)$ and then for regular hyperbolic elements in $\PSp(n,1)$. In
Section~\ref{sp}, we study the $\PSp(1,1)$ case when the two elements have common
fixed points. Finally, Section~\ref{qnt} gives a complete explicit matrix criterion
for strong double reversibility in $\PSp(1,1)$.

\section{Preliminaries} \label{pre}

\subsection{Doubly reversible pairs}
Let $G$ act on $G\times G$ by simultaneous conjugation. For $(g_1,g_2)\in G^2$ set
\[
S_G(g_1,g_2)=Z_G(g_1)\cap Z_G(g_2)
\]
and
\[
R_G(g_1,g_2)=\left\{h\in G:
\begin{array}{l}
hg_1h^{-1}=g_1^{-1},\\[-2pt]
hg_2h^{-1}=g_2^{-1}
\end{array}\right\}.
\]
Define $E_G(g_1,g_2)=S_G(g_1,g_2)\cup R_G(g_1,g_2)$.

An element of $R_G(g_1,g_2)$ is called a \emph{reverser} of $(g_1,g_2)$.  Thus $(g_1,g_2)$ is doubly reversible precisely when
$R_G(g_1,g_2)\neq\emptyset$, and strongly doubly reversible precisely when
$R_G(g_1,g_2)$ contains an involution.
\begin{lemma}\label{drl}
$E_G(g_1,g_2)$ is a subgroup of $G$. Moreover,
$S_G(g_1,g_2)$ is a normal subgroup of $E_G(g_1,g_2)$ of index at most two.
\end{lemma}
\begin{proof}
Products of two reversers lie in $S_G(g_1,g_2)$, while multiplying a stabilizer and a
reverser gives a reverser. Thus $E_G(g_1,g_2)$ is a subgroup. If
$R_G(g_1,g_2)=\emptyset$, then $E_G=S_G$. Otherwise define
$\phi:E_G\to\{\pm1\}$ by $\phi(h)=1$ on $S_G$ and $\phi(h)=-1$ on $R_G$.
This is a surjective homomorphism with kernel $S_G$, proving the assertion.
\end{proof}

The following elementary observation will also be useful. It records precisely how a
common reverser may be modified by the common centralizer.
\begin{proposition}\label{prop:coset}
Assume $R_G(g_1,g_2)\neq\emptyset$ and that $g_1\neq g_1^{-1}$ or
$g_2\neq g_2^{-1}$. Fix $C\in R_G(g_1,g_2)$ and put $Z=S_G(g_1,g_2)$. Then
$C^2\in Z$, $C$ normalizes $Z$, and
\[
R_G(g_1,g_2)=ZC=CZ.
\]
Consequently $(g_1,g_2)$ is strongly doubly reversible if and only if there is
$z\in Z$ such that
\[
z(CzC^{-1})=C^{-2}.
\]
\end{proposition}
\begin{proof}
Conjugating twice by $C$ shows $C^2\in Z$, and a direct calculation shows that
$C Z C^{-1}=Z$. If $D$ is another common reverser, then $C^{-1}D\in Z$, so
$R_G=ZC=CZ$. Finally,
\[
(zC)^2=z(CzC^{-1})C^2,
\]
which gives the last assertion.
\end{proof}

In particular, if the common centralizer is trivial, a common reverser is unique and is
automatically an involution.

\subsection{The Quaternions}
Let $\mathbb{H} := \mathbb{R} + \mathbb{R}i + \mathbb{R}j + \mathbb{R}k$ denote the division algebra of Hamilton’s quaternions, where the fundamental relations are given by 
$i^2 = j^2 = k^2 = ijk = -1.$ 
Every element of $\mathbb{H}$ can be written uniquely in the form 
$q = a + bi + cj + dk, ~\textit{where}~ a,b,c,d \in \mathbb{R}.$  
Alternatively, viewing $\mathbb{H}$ as a two-dimensional vector space over $\mathbb{C}$, we may express $q = c_1 + c_2 j, ~\textit{with} ~c_1, c_2 \in \mathbb{C}.$ 
The modulus (or norm) of $q$ is defined by
$|q| = \sqrt{a^2 + b^2 + c^2 + d^2}.$
We denote the set
$$\Sp(1) := \{ q \in \mathbb{H} : |q| = 1 \}$$
by the group of unit quaternions.

\medskip

We consider $\H^n$ as a right vector space over the quaternions. A non-zero vector $v \in \H^n $ is said to be a (right) eigenvector of $A$ corresponding to a (right) eigenvalue $\lambda \in \H $ if the equality $ A v = v\lambda $ holds.

Eigenvalues of every matrix over the quaternions occur in similarity classes, and each similarity class of eigenvalues contains a unique complex number with non-negative imaginary part. Here, instead of similarity classes of eigenvalues, we will consider the \textit{unique complex representatives} with non-negative imaginary parts as eigenvalues unless specified otherwise. In places where we need to distinguish between the similarity class and a representative, we shall write the similarity class of an eigenvalue representative $\lambda$ as $[\lambda]$. 

\subsection{Quaternionic Hyperbolic Space}
Let $\V=\H^{n,1}$ denote the right vector space of dimension $n+1$ over $\H$, equipped with the Hermitian form:  
$$\langle\z,\w\rangle=\w^{\ast}H_1\z=\bar w_{n+1} z_1+\Sigma_{i=2}^n \bar w_i z_i+\bar w_1 z_{n+1}, $$
where $\ast$ denotes the conjugate transpose, and 
\begin{center}
$$H_1= \begin{pmatrix}
            0 & 0 & 1\\
           0 & I_{n-1} & 0 \\
 1 & 0 & 0
          \end{pmatrix} .$$
\end{center}
We consider the following subspaces of $\H^{n,1}:$
$$\V_{-}=\{\z\in\H^{n,1}:\langle\z,\z \rangle<0\}, ~ \V_+=\{\z\in\H^{n,1}:\langle\z,\z \rangle>0\},$$
$$\V_{0}=\{\z \in\H^{n,1} \setminus \{{ 0}\}:\langle\z,\z \rangle=0\}.$$
Let $\P:\H^{n,1} \setminus \{0\}\longrightarrow  \H \P^n$ be the right projection onto the quaternionic projective space. The image of a vector $\z$ will be denoted by $z$.   

The projective model of the quaternionic hyperbolic space is given by $\h^n=\P(\V_-)$. The boundary at infinity of this space is $\partial\h^n=\P(\V_0)$.

The above Hermitian form may be replaced by an equivalent one associated with the $(n+1) \times (n+1)$ matrix $H_o$: 
$$H_o={\rm diag }~(-1, 1, \ldots, 1). $$
where the corresponding Hermitian form $\langle \z, \w \rangle_o=\w^{\ast} H_o \z$ gives the ball model of $\h^n$. 

Given a point $z$ of $\h^n \cup \partial \h^n -\{\infty\} \subset\H \P^n$ we may lift $z=(z_1,z_2,\ldots,z_n)$ to a point $\z$ in $\V_0 \cup \V_{-}$, given by 
 $$\z=\begin{pmatrix}
                z_1\\z_2\\ \vdots\\1
               \end{pmatrix}.$$
Here $\z$ is called the \emph{standard lift} of $z$.  There are two points: `zero' and `infinity'  in the boundary  given by:
$$o=\begin{pmatrix}
                0\\0\\\vdots \\1
               \end{pmatrix}, ~~~~ \infty=\begin{pmatrix}
                1\\0\\\vdots\\0
               \end{pmatrix}.$$

Let ${\rm Sp}(n,1)$ be the isometry group of the Hermitian form $H_1$. Each matrix $A$ in ${\rm Sp}(n,1)$ satisfies the relation $A^{-1}=
{H_1}^{-1}A^{\ast}H_1$, where $A^{\ast}$ is the conjugate transpose of $A$. The isometry group of  $\h^n$ is the projective unitary group ${\rm PSp}(n,1)={\rm Sp }(n,1)/\{\pm I\}$. However, we shall mostly deal with the linear group $\Sp(n,1)$. 


\subsection{Classification of elements} Following the terminology in \cite{Rat}, recall that an element $g \in \Sp(n,1)$ is called \emph{hyperbolic} if it has exactly two fixed points on the boundary. 

 An element $g \in \Sp(n,1)$ is called parabolic if it has a unique fixed point on the boundary, and elliptic if it has a fixed point in $\h^n$. An element $g$ in $\Sp(n,1)$ belongs to exactly one of these three classes. 

\subsection{ Hyperbolic Isometries }   \label{ltr}
Consider a hyperbolic isometry $A \in \mathrm{Sp}(n,1)$. Let $[\lambda]$ denote the conjugacy class of eigenvalues associated with $A$, and choose a representative eigenvalue $\lambda$ with a corresponding eigenvector $\x$. The vector $\x$ determines a point in quaternionic projective space $\mathbb{H} \mathrm{P}^n$, which lies either on the boundary $\partial \h^n$ or is a point in $\P(\V_{+})$. The corresponding line $\x\mathbb{H}$ in the space $\mathbb{H}^{n,1}$ represents the lift of this projective point and is invariant under the action of $A$. This line is the eigenspace generated by $\x$.

In the hyperbolic case, two of the eigenvalue classes are of null-type, with their associated eigenlines corresponding to fixed points on the boundary: one attracting, the other repelling. Suppose the repelling fixed point of $A$ on $\partial \mathbb{H}^n$ is denoted by $r_A$ and corresponds to the eigenvalue $re^{i\theta}$, while the attracting fixed point $a_A$ corresponds to the eigenvalue $r^{-1}e^{i\theta}$. Let $\mathbf{r}_A$ and $\mathbf{a}_A$ denote their respective lifts to $\mathbb{H}^{n,1}$. Additionally, for each $j$, let $\x_{j,A}$ be an eigenvector of $A$ associated with the eigenvalue $e^{i\phi_j}$. It is convenient to assume that the angles $\theta, \phi_j$ lie within the interval $[0, \pi]$ and $0<r<1$. Each $\x_{j,A}$ defines a point in $\mathbb{P}(V_+)$.

Given $(r,\theta,\phi_1,\ldots,\phi_{n-1})$, let $E_A$ denote the diagonal matrix, with respect to the standard Hermitian form $H_1$,
\begin{equation}
E_A=\operatorname{Diag}\bigl(re^{i\theta},e^{i\phi_1},\ldots,e^{i\phi_{n-1}},r^{-1}e^{i\theta}\bigr).
\end{equation}

Construct the matrix
$$
C_A = \left[ {\r}_A \quad \x_{1,A} \quad \cdots \quad \x_{n-2,A} \quad \x_{n-1,A} \quad {\a_A}\right],
$$
whose columns are the eigenvectors corresponding to the eigenvalues used in $E_A$. By suitably scaling the eigenvectors, we can ensure that $C_A$ belongs to $\mathrm{Sp}(n,1)$, by enforcing the normalization: 
$$
\langle {\a}_A, {\r}_A \rangle = 1, \quad \langle \x_{j,A}, \x_{j,A} \rangle = 1.
$$
With this choice of basis, the matrix $A$ is conjugate to the diagonal matrix $E_A$, i.e.,
$$
A = C_A E_A C_A^{-1}.
$$
So, every hyperbolic element  $A$ in $\Sp(n,1)$ is conjugate to a matrix of the form  $E_A$.
\begin{lemma}[Chen-Greenberg] \cite{chen}\label{hycl}
 Two hyperbolic elements in ${\rm Sp}(n,1)$ are conjugate if and only if they have the same similarity classes of eigenvalues.
\end{lemma}

\subsection{Reversers of hyperbolic elements}
We shall repeatedly use the following elementary fact.
\begin{lemma}\label{loxo1}
Let $A\in\PSp(n,1)$ be hyperbolic. If $C\in\PSp(n,1)$ satisfies
$CAC^{-1}=A^{-1}$, then $C$ interchanges the attracting and repelling fixed points of
$A$. In particular, if a lift of $A$ is written as a product of two skew-involutions,
then each of those skew-involutions interchanges the two fixed points.
\end{lemma}
\begin{proof}
The attracting fixed point of $A$ is carried by $C$ to the attracting fixed point of
$A^{-1}$, which is the repelling fixed point of $A$, and conversely. If
$A=i_1i_2$ with $i_1^2=i_2^2=-I$, then $i_jAi_j^{-1}=A^{-1}$ for $j=1,2$, so the
first assertion applies.
\end{proof}

\section{Strongly doubly reversible pairs in $\PSp(1)$}\label{sp11} 
\begin{lemma}\label{lem}
	Let $q \in \Sp(1)$ be such that $e^{-i\theta}= q e^{i \theta} q^{-1}$ where $\theta \neq 0, \pi$ then $q= e^{i\phi}j$ for some $\phi \in [0, 2\pi).$
\end{lemma}
\begin{proof}
	Let $q \in \Sp(1)$ such that $q = c_1+c_2j$ for $c_1,c_2 \in \C.$
	We have: $$qe^{i\theta} = (c_1+c_2j) e^{i\theta} = c_1 e^{i \theta} + c_2e^{-i\theta}j,$$
	$$e^{-i\theta} q = e^{-i\theta}(c_1+c_2j) = c_1 e^{-i \theta} + c_2e^{-i\theta}j.$$
	Comparing the two sides, we get $c_1 e^{i \theta} = c_1 e^{-i \theta}$. Since $\theta \neq 0, \pi$, this implies $c_1 = 0$. As $q \in \Sp(1),$ we then have $q= e^{i\phi}j$ for some $\phi \in [0, 2\pi)$. 
\end{proof}
\begin{remark}
	In the above lemma,  $q= e^{i\phi}j$ for $\phi \in [0, 2\pi),$ and hence $q^2= -1$. Thus, $q$ is an involution in $\PSp(1).$ Consequently, every element in $\PSp(1)$ is strongly reversible.
\end{remark}

\begin{theorem}\label{sp1} 
Any two elements in $\PSp(1)$ are strongly doubly reversible.
\end{theorem}
\begin{proof} Let $p_1$ and $p_2$ be elements in $\Sp(1)$. 
	Without loss of generality, assume $p_1= e^{i\theta}$ and $p_2= c_1+c_2j$ be elements in $\Sp(1)$. Then we need to find $q$ such that  ${p_1}^{-1}= e^{-i\theta}= q e^{i \theta} q^{-1}, ~p_2^{-1}=q p_2 q^{-1},$ where $q^2=\pm 1$. 
	
	Now observe that, by using lemma \ref{lem}, ${p_1}^{-1}= e^{-i\theta}= q e^{i \theta} q^{-1},$ holds for any $q= e^{i\theta_1}j$ where $\theta_1 \in [0, 2\pi)$. So, we have $$qp_2= e^{i \theta_1}j(c_1+c_2j)= e^{i \theta_1}jc_1 + e^{i \theta_1}j c_2j= e^{i \theta_1}\bar c_1j- e^{i \theta_1}\bar c_2,$$
	$$p_2^{-1} q = (\bar c_1- c_2j)e^{i \theta_1}j= \bar c_1 e^{i \theta_1}j- c_2j e^{i \theta_1}j= e^{i \theta_1}\bar c_1j+ c_2e^{-i \theta_1}.$$ For $p_1$ and $p_2$ strongly doubly reversible by $q$, we require $c_2 e^{-i \theta_1}= - e^{i \theta_1}\bar c_2,$ which is equivalent to $\operatorname{Re}(c_2 e^{-i \theta_1})=0$. 
    
    Writing $c_2=c+di,$ this gives $\cos{\theta_1} ~c+\sin{\theta_1}~ d =0,$ which always has a solution $\theta_1 \in [0, 2\pi).$ That means if we know $c_2$, we can always find $\theta_1$ such that $\operatorname{Re}(c_2 e^{-i \theta_1})=0.$ 
	Then we get $${p_1^{-1}}= e^{-i\theta}= q e^{i \theta} q^{-1}, ~p_2^{-1}=q p_2 q^{-1},$$ where $q= e^{i \theta_1}j$ for $\theta_1 \in [0, 2\pi)$. That means, $p_1$ strongly doubly reversible with $p_2$ in $\PSp(1)$ via $q$.
\end{proof}

\begin{corollary}
	Every pair of elements in ${\rm SO}(3)$ is strongly doubly reversible. 
\end{corollary}
\begin{proof}
It is a well-known result that $\Sp(1)$ is a double cover of ${\rm SO}(3)$. Hence, ${\rm PSp}(1) \simeq {\rm SO}(3)$, and the result follows from Theorem~\ref{sp1}.
\end{proof}

As an application of the above theorem, we provide a simple proof of the following result (see in \cite{bm}).

\begin{theorem}
	Every pair of elements in ${\rm SO}(4)$ is strongly doubly 
    reversible. 
\end{theorem}
\begin{proof}
It is a well-known result that $\Sp(1) \times \Sp(1)$ is a double 
cover of ${\rm SO}(4)$. Let $A, B \in {\rm SO}(4)$, and let 
$\tilde{A}= (A_1,A_2)$, $\tilde{B}= (B_1,B_2)$ be their respective 
lifts in $\Sp(1) \times \Sp(1)$. By Theorem~\ref{sp1}, for each 
$i=1,2,$ there exist exist skew-involutions $\tilde{\alpha}_i, 
\tilde{\beta}_i, \tilde{\gamma}_i$ such that 
$$A_i= \tilde{\alpha_i} \tilde{\beta_i}, \quad B_i=\tilde{\beta_i}  
\tilde{\gamma_i} .$$ 
Define $\tilde{\alpha}= (\tilde{\alpha_1},\tilde{\alpha_2})$ and 
$\tilde{\beta}= (\tilde{\beta_1},\tilde{\beta_2})$ then we get 
$\tilde{\alpha} \tilde{\beta}= 
(\tilde{\alpha_1}\tilde{\beta_1},\tilde{\alpha_2} \tilde{\beta_2}) =
(A_1,A_2)=\tilde{A}.$ Taking projection, we obtain $$ A= 
\pi(\tilde{\alpha} \tilde{\beta})=\pi(\tilde{\alpha}) 
\pi(\tilde{\beta}) =\alpha \beta,$$ where $\alpha= 
\pi(\tilde{\alpha}), \beta= \pi(\tilde{\beta}) \in {\rm SO}(4).$ We 
can observe that $$\alpha^2= \pi(\tilde{\alpha})\pi(\tilde{\alpha}) = 
\pi(\tilde{\alpha}^2)= \pi((-1,-1))=1,$$ so $\alpha$ is an involution, 
and similarly, $\beta$ is an involution. By the same argument, 
$B=\beta \gamma$ with $\beta^2=\gamma^2= 1$. Thus, $(A,B)$ is strongly 
doubly reversible.
\end{proof}

\section{Strongly Doubly Reversible Pairs in $\PSp(n,1)$}\label{spn1s}
Passing from $\Sp(n,1)$ to $\PSp(n,1)=\Sp(n,1)/\{\pm I\}$ does not affect
Haar measure-zero statements, since the quotient is by the finite central
subgroup $\{\pm I\}$.
\subsection{Proof of Theorem~\ref{spn1}}

Let \(G=\PSp(n,1)\) and \(\mathfrak g=\mathfrak{sp}(n,1)\). Recall that
$$
\dim_{\mathbb R}G=(n+1)(2n+3).
$$
Let \(s\in G\) be an involution. A lift \(\widetilde s\in\Sp(n,1)\) satisfies
\(\widetilde s^{\,2}=\pm I\). Since \(\Ad(s)^2=I\), we have
$$
\mathfrak g=\mathfrak g_{+1}(s)\oplus\mathfrak g_{-1}(s),
$$
where
$$
\mathfrak g_{\pm1}(s)=\{X\in\mathfrak g:\Ad(s)X=\pm X\}.
$$
The \(+1\)-eigenspace is the Lie algebra of \(Z_G(s)\).

If \(\widetilde s^{\,2}=-I\), then
\(\mathfrak g_{+1}(s)\cong\mathfrak u(n,1)\), and hence
$$
\dim_{\mathbb R}\mathfrak g_{-1}(s)
=(n+1)(2n+3)-(n+1)^2
=(n+1)(n+2).
$$
If \(\widetilde s^{\,2}=I\), the \(+1\)- and \(-1\)-eigenspaces of
\(\widetilde s\) are orthogonal and nondegenerate. If their quaternionic
dimensions are \(p\) and \(q\), with \(p+q=n+1\), then
$$
\dim_{\mathbb R}\mathfrak g_{-1}(s)
=4pq\le (n+1)^2<(n+1)(n+2).
$$
Thus, for every involution \(s\in G\),
\begin{equation}\label{eq:max-minus}
\dim_{\mathbb R}\mathfrak g_{-1}(s)\le (n+1)(n+2).
\end{equation}

The conjugacy class of \(s\) has dimension
$$
\dim G-\dim Z_G(s)=\dim\mathfrak g_{-1}(s).
$$
There are only finitely many conjugacy types of involutions in \(G\):
when \(\widetilde s^{\,2}=I\), the type is determined by the dimensions
and signatures of its two eigenspaces, while the case
\(\widetilde s^{\,2}=-I\) gives the complex-structure type described above.
Consequently, if
$$
\mathcal I_0=\{s\in G:s^2=1\}, \hbox{ then we have} 
$$
\begin{equation}\label{eq:inv-dim}
\dim\mathcal I_0\le (n+1)(n+2).
\end{equation}
Note that \(\mathcal I_0\) consists of the involutions together with the identity.

Fix \(h\in\mathcal I_0\) and set
$$
R_h=\{g\in G:hgh^{-1}=g^{-1}\}.
$$
Since \(h^2=1\),
$$
hgh^{-1}=g^{-1}
\quad\Longleftrightarrow\quad
(hg)^2=1.
$$
Therefore, 
\begin{equation}\label{eq:Rh-inv}
R_h=h\mathcal I_0.
\end{equation}
Left multiplication by \(h\) is a diffeomorphism of \(G\), so
$$
\dim R_h=\dim\mathcal I_0\le (n+1)(n+2).
$$
We now allow \(h\) to vary. Consider
$$
\mathcal R=
\{(g_1,g_2,h)\in G^2\times\mathcal I_0:
hg_ih^{-1}=g_i^{-1},\ i=1,2\}.
$$
For each \(h\in\mathcal I_0\), the fiber in the
\((g_1,g_2)\)-variables is \(R_h\times R_h\), and hence has dimension at
most \(2(n+1)(n+2)\). By \eqref{eq:inv-dim},
$$
\dim\mathcal R\le 3(n+1)(n+2).
$$
On the other hand,
$$
\dim(G\times G)=2(n+1)(2n+3),
$$
and hence, for every \(n\ge1\),
$$
\begin{aligned}
\dim(G\times G)-\dim\mathcal R
&\ge
2(n+1)(2n+3)-3(n+1)(n+2)\\
&=n(n+1)>0.
\end{aligned}
$$
To justify the dimension argument, we may realize \(G\) in its adjoint
representation as a real algebraic group. Then \(\mathcal I_0\),
\(\mathcal R\), and the equations defining the reversing relations are
semialgebraic. The strongly doubly reversible locus is contained in the
projection of \(\mathcal R\) onto the first two factors, and semialgebraic
projection does not increase dimension. Therefore this locus has
dimension strictly smaller than \(\dim(G\times G)\), and hence has Haar
measure zero. \qed
\subsection{Proof of Corollary~1.3}

Let \(G=\SU(n,1)\), \(n\ge2\), and put \(N=n+1\). Then
\(\dim_{\mathbb R}G=N^2-1\). If \(s\in G\) is an involution whose
\((-1)\)-eigenspace has complex dimension \(m\), then its conjugacy class has
real dimension \(2m(N-m)\). Hence
$$
\dim\mathfrak g_{-1}(s)\le
\left\lfloor\frac{N^2}{2}\right\rfloor
=
\left\lfloor\frac{(n+1)^2}{2}\right\rfloor.
$$
There are only finitely many conjugacy types of involutions. Thus, for
\(\mathcal I_0=\{s\in G:s^2=1\}\),
\[
\dim\mathcal I_0\le
\left\lfloor\frac{(n+1)^2}{2}\right\rfloor.
\]
For fixed \(h\in\mathcal I_0\), we have
\(R_h=\{g\in G:hgh^{-1}=g^{-1}\}=h\mathcal I_0\), and therefore
\[
\dim R_h\le
\left\lfloor\frac{(n+1)^2}{2}\right\rfloor.
\]
Consequently, the corresponding incidence set of strongly doubly reversible
pairs has dimension at most
\(3\left\lfloor (n+1)^2/2\right\rfloor\). Since, for \(n\ge2\),
$$
3\left\lfloor\frac{(n+1)^2}{2}\right\rfloor
<2\bigl((n+1)^2-1\bigr)=2\dim\SU(n,1),
$$
the strongly doubly reversible pairs form a Haar measure zero subset of
\(\SU(n,1)\times\SU(n,1)\). \qed

\subsection{Proof of Corollary~1.4}

Let \(G=\PSp(n)\), \(n\ge2\). Then
\(\dim_{\mathbb R}G=n(2n+1)\). If an involution \(s\in G\) has a lift
satisfying \(\widetilde s^{\,2}=-I\), then
\(\mathfrak g_{+1}(s)\cong\mathfrak u(n)\), and hence
\[
\dim\mathfrak g_{-1}(s)=n(2n+1)-n^2=n(n+1).
\]
If \(\widetilde s^{\,2}=I\), with quaternionic eigenspace dimensions \(p\)
and \(n-p\), then
\[
\dim\mathfrak g_{-1}(s)=4p(n-p)\le n^2<n(n+1).
\]
Thus every involution conjugacy class has dimension at most \(n(n+1)\).
Since there are only finitely many conjugacy types,
\(\dim\mathcal I_0\le n(n+1)\), where
\(\mathcal I_0=\{s\in G:s^2=1\}\). For every \(h\in\mathcal I_0\),
\(R_h=h\mathcal I_0\), so \(\dim R_h\le n(n+1)\). Hence the incidence set
of strongly doubly reversible pairs has dimension at most \(3n(n+1)\).
Since, for \(n\ge2\),
$$
3n(n+1)<2n(2n+1)=2\dim G,
$$
the strongly doubly reversible pairs form a Haar measure zero subset of
\(\PSp(n)\times\PSp(n)\). \qed
\subsection{Proof of Corollary~ \ref{cor:SO}}
Let $G=\SO_0(n,1)$, so that $\dim G=n(n+1)/2$. Every involution
$s\in G$ has eigenvalues $\pm1$. If $m=\dim E_-(s)$, then
\[
\dim Z_G(s)
=\frac{m(m-1)}2+\frac{(n+1-m)(n-m)}2,
\]
and hence its conjugacy class has dimension $m(n+1-m)\le
\lfloor (n+1)^2/4\rfloor$. Since there are only finitely many conjugacy
types of involutions, $\dim I_0\le\lfloor (n+1)^2/4\rfloor$.

For fixed $h\in I_0$, the set of elements reversed by $h$ is $R_h=hI_0$, and
therefore $\dim R_h\le\lfloor (n+1)^2/4\rfloor$. Thus the incidence set
\[
R=\{(g_1,g_2,h)\in G^2\times I_0:
hg_i h^{-1}=g_i^{-1},\ i=1,2\}
\]
satisfies $\dim R\le3\lfloor (n+1)^2/4\rfloor$. Since $n\ge4$,
\[
3\left\lfloor\frac{(n+1)^2}{4}\right\rfloor<n(n+1)
=\dim(G\times G).
\]
Hence the projection of $R$ onto $G\times G$, and therefore the strongly
doubly reversible locus, has Haar measure zero.

The compact case $G=\SO(n+1)$ is identical, since
$\dim\SO(n+1)=n(n+1)/2$ and the involutions have the same conjugacy-class
dimension bound.\qed
\subsection{Strongly \(k\)-reversible tuples}

The same argument applies to \(k\)-tuples. Let
\[
\mathcal R_k=\{(g_1,\ldots,g_k,h)\in G^k\times\mathcal I_0:
hg_ih^{-1}=g_i^{-1},\ 1\le i\le k\}.
\]
For fixed \(h\), the fiber is \(R_h^k\), and hence has dimension at most
\(k(n+1)(n+2)\). Since
\(\dim\mathcal I_0\le(n+1)(n+2)\),
\[
\dim\mathcal R_k\le(k+1)(n+1)(n+2).
\]
Moreover,
$$
\begin{aligned}
k\dim G-\dim\mathcal R_k
&\ge k(n+1)(2n+3)-(k+1)(n+1)(n+2)\\
&=(n+1)\bigl[n(k-1)+(k-2)\bigr]>0
\end{aligned}
$$
for every \(n\ge1\) and \(k\ge2\). Therefore the projection of
\(\mathcal R_k\) onto \(G^k\) has dimension strictly smaller than
\(\dim G^k\), and hence has Haar measure zero.

\begin{theorem}\label{thm:krev}
For every \(n\ge1\) and \(k\ge2\), the set of strongly \(k\)-reversible
tuples in \(\PSp(n,1)^k\) has Haar measure zero.
\end{theorem}
\begin{remark}\label{rk}
The compact group $\PSp(1)\cong{\rm SO}(3)$ is exceptional: by
Section~\ref{sp11}, every pair in $\PSp(1)$ is strongly doubly reversible. The
corresponding compact-group dimension inequality becomes an equality when $n=1$.
\end{remark}

\section{Double versus strong double reversibility for hyperbolic pairs}\label{sp-1}

We begin with the elementary $2\times2$ calculation that underlies both the
one-dimensional and the regular higher-dimensional arguments.
\begin{lemma}\label{invv}
Let

$$
A=\begin{pmatrix}re^{i\theta}&0\\0&r^{-1}e^{i\theta}\end{pmatrix}\in\Sp(1,1),
\qquad 0<r<1,\quad \theta\in[0,\pi],
$$

and let $C\in\Sp(1,1)$ satisfy $CAC^{-1}=\varepsilon A^{-1}$, where
$\varepsilon\in \{\pm1\}$. Then

$$
C=C_y:=\begin{pmatrix}0&y\\ \bar y^{-1}&0\end{pmatrix},
\qquad y\in\H^*,
$$

and, writing $u=y/|y|$, we have $C_y^2=u^2I$. Moreover,

$$
ye^{i\theta}=\varepsilon e^{-i\theta}y,
$$

and:
\begin{enumerate}
\item[(i)] if $\varepsilon=1$ and $\theta\in(0,\pi)$, then
$y\in\C^*j$ and $C_y^2=-I$;
\item[(ii)] if $\varepsilon=1$ and $\theta\in (0,\pi)$, there is no
restriction on $y\in\H^*$;
\item[(iii)] if $\varepsilon=-1$, then necessarily $\theta=\pi/2$ and
$y\in\C^*$.
\end{enumerate}
\end{lemma}
\begin{proof}
First suppose $\varepsilon=1$. Write
$C=\left(\begin{smallmatrix}x&y\\ z&w\end{smallmatrix}\right)$. From
$CA=A^{-1}C$ we obtain
$$
x=r^2e^{i\theta}xe^{i\theta},\qquad
w=r^{-2}e^{i\theta}we^{i\theta}.
$$
Taking norms gives $x=w=0$. Since $C\in\Sp(1,1)$, we have
$z=\bar y^{-1}$, so
$$
C=C_y=\begin{pmatrix}0&y\\ \bar y^{-1}&0\end{pmatrix},
\qquad
ye^{i\theta}=e^{-i\theta}y.
$$
Writing $y=|y|u$, with $u\in\Sp(1)$, gives
$$
C_y^2=\begin{pmatrix}
y\bar y^{-1}&0\\0&\bar y^{-1}y
\end{pmatrix}=u^2I.
$$
If $\theta\in(0,\pi)$, writing $y=c_1+c_2j$, $c_1,c_2\in\C$, the above
relation forces $c_1=0$. Thus $y\in\C^*j$, so $u^2=-1$ and hence
$C_y^2=-I$. If $\theta=0$ or $\pi$, the relation imposes no restriction
on $y$.

Now suppose $\varepsilon=-1$. The same norm argument gives $x=w=0$ and
hence $C=C_y$, while the remaining equation is
$$
ye^{i\theta}=-e^{-i\theta}y.
$$
Writing again $y=c_1+c_2j$ gives $c_2=0$ and
$$
c_1e^{i\theta}=-e^{-i\theta}c_1.
$$
Since $y\ne0$, we have $c_1\ne0$, and therefore
$e^{2i\theta}=-1$. Thus $\theta=\pi/2$ and $y\in\C^*$.
\end{proof}

\subsection{Proof of Theorem~\ref{hyp}}
Together, Theorem~\ref{thm:n1equiv} and Proposition~\ref{prop:regular} will prove
Theorem~\ref{hyp}.

\begin{theorem}\label{thm:n1equiv}
Let $A \in \PSp(1,1)$ be hyperbolic and let $B \in \PSp(1,1)$ be arbitrary.
Then $(A,B)$ is doubly reversible if and only if it is strongly doubly reversible.
\end{theorem}
\begin{proof}
Only the converse requires proof. Suppose $(A,B)$ is doubly reversible.
Conjugating in $\Sp(1,1)$ and choosing lifts, we may assume
\[
A=\operatorname{diag}(re^{i\theta},r^{-1}e^{i\theta}),
\qquad 0<r<1,\quad \theta\in[0,\pi],
\]
and that a common reverser $C\in\Sp(1,1)$ satisfies
\[
CAC^{-1}=\varepsilon A^{-1},\qquad
CBC^{-1}=\delta B^{-1},
\]
with $\varepsilon,\delta\in\{\pm1\}$. By Lemma~\ref{invv},
$C=C_y$ and $C^2=u^2I$, where $u=y/|y|$.

Since $A$ is hyperbolic, $A\neq A^{-1}$, so Proposition~\ref{prop:coset}
applies to the pair $(A,B)$ in $\PSp(1,1)$: writing $Z=S_{\PSp(1,1)}(A,B)$, it
suffices to produce $z\in Z$ with
\[
z\,(CzC^{-1})=C^{-2}.
\]

If $u^2=\pm1$, then $C^2=\pm I$, so $C^{-2}=1$ in $\PSp(1,1)$ and $z=1$ has the
required property. This includes, in particular, the case $\varepsilon=1$ and
$\theta\in(0,\pi)$, where Lemma~\ref{invv}(i) gives $u^2=-1$.

Assume therefore that $u^2\notin\R$. Then $u\notin\R$. By
Lemma~\ref{invv}, we are in one of the following two cases:
either $\varepsilon=1$ and $\theta\in\{0,\pi\}$, or
$\varepsilon=-1$, $\theta=\pi/2$, and $y\in\C^\ast$.

Conjugating the relation $CBC^{-1}=\delta B^{-1}$ once more gives
$C^2BC^{-2}=B$. Since $C^2=u^2I$, this means that every entry of $B$
commutes with $u^2$. As $u^2\notin\R$, we have
$Z_{\H}(u^2)=Z_{\H}(u)$, and hence every entry of $B$ commutes with $u$.

Put $z=u^{-1}I$. Since
\[
H_1=\begin{pmatrix}0&1\\1&0\end{pmatrix}
\]
has real entries and $|u|=1$, we have $z^*H_1z=H_1$, and hence
$z\in\Sp(1,1)$. Moreover, $z$ commutes with $B$. It also commutes with
$A$: in the case $\varepsilon=1$, $\theta\in\{0,\pi\}$, the entries of
$A$ are real, while in the case $\varepsilon=-1$, $\theta=\pi/2$, the
entries of $A$ are $ri$ and $r^{-1}i$, and $u\in\C$. Thus $z$ projects into
$Z$.

Finally, $z$ commutes with $C$. Indeed, writing $y=|y|u$, we have
$u^{-1}y=yu^{-1}=|y|$, while
$\bar y^{-1}=|y|^{-1}u$, and hence
\[
u^{-1}\bar y^{-1}=\bar y^{-1}u^{-1}=|y|^{-1}.
\]
Therefore $CzC^{-1}=z$, and
\[
z\,(CzC^{-1})=z^2=u^{-2}I=(u^2I)^{-1}=C^{-2}.
\]
By Proposition~\ref{prop:coset}, $(A,B)$ is strongly doubly reversible.
\end{proof}
\subsection{Regular hyperbolic elements in $\PSp(n,1)$}
A hyperbolic element $A\in\PSp(n,1)$ has a lift conjugate to
\[
D_A=\operatorname{diag}\bigl(re^{i\theta},e^{i\phi_1},\ldots,e^{i\phi_{n-1}},
r^{-1}e^{i\theta}\bigr),\qquad 0<r<1.
\]
\begin{definition}\label{def:regular}
We call $A$ \emph{regular} here if
\[
\theta\in(0,\pi)\setminus\{\pi/2\},\qquad
\phi_j\in(0,\pi),
\]
and the similarity classes $[e^{i\phi_1}],\ldots,[e^{i\phi_{n-1}}]$ are pairwise
distinct.
\end{definition}

\begin{proposition}\label{prop:regular}
If $A \in \PSp(n,1)$ is regular hyperbolic, then every projective reverser of $A$
is an involution. Consequently, for any $B \in \PSp(n,1)$, if $(A,B)$ is doubly
reversible then it is strongly doubly reversible.
\end{proposition}
\begin{proof}
Let $C$ reverse $A$ and choose lifts $\widetilde A,\widetilde C\in\Sp(n,1)$ such that
\[
\widetilde C\widetilde A\widetilde C^{-1}=\varepsilon\widetilde A^{-1},
\qquad \varepsilon\in\{\pm1\}.
\]
If $\varepsilon=-1$, then $\widetilde A$ and $-\widetilde A^{-1}$ are conjugate in
$\Sp(n,1)$, hence have the same eigenvalue classes. Recall that a quaternionic
similarity class $[\lambda]$ is determined by $|\lambda|$ and $\Re\lambda$. The
classes of $\widetilde A$ are $[re^{i\theta}]$, $[r^{-1}e^{i\theta}]$ and
$[e^{i\phi_j}]$, of moduli $r$, $r^{-1}$ and $1$; those of $-\widetilde A^{-1}$
are $[-re^{-i\theta}]$, $[-r^{-1}e^{-i\theta}]$ and $[-e^{-i\phi_j}]$, of the
same moduli. Since $0<r<1$, each of $\widetilde A$ and $-\widetilde A^{-1}$ has
exactly one eigenvalue class of modulus $r$, so
$[re^{i\theta}]=[-re^{-i\theta}]$. Comparing real parts gives
$r\cos\theta=-r\cos\theta$, hence $\theta=\pi/2$, contrary to regularity. Thus
$\varepsilon=1$.

Conjugate so that $\widetilde A=D_A$. The reverser interchanges the two null
eigenlines and, because the unit eigenvalue classes are pairwise distinct, preserves
each positive eigenline. On the null $2\times2$ block Lemma~\ref{invv} gives square
$-I$. On a positive eigenline with eigenvalue $e^{i\phi_j}$, write
$\widetilde Cx_j=x_jq_j$. The relation
$\widetilde C\widetilde A=\widetilde A^{-1}\widetilde C$ gives
\[
q_je^{i\phi_j}=e^{-i\phi_j}q_j.
\]
Since $|q_j|=1$ and $\phi_j\in(0,\pi)$, we have $q_j\in\C j$ and $q_j^2=-1$.
Hence $\widetilde C^2=-I$ on every eigenspace, so $C^2=1$ in $\PSp(n,1)$.
A common reverser of $(A,B)$ in particular reverses the regular member and therefore is
an involution.
\end{proof}

\section{Strongly doubly reversible hyperbolic pairs in $\PSp(1,1)$}\label{sp}
\subsection{Pairs with a common fixed point}
\begin{prop}\label{common point}
Let $A$ and $B$ be hyperbolic elements in $\Sp(1,1)$ with a common fixed point. Then
$A$ and $B$ are strongly doubly reversible in $\PSp(1,1)$ if and only if their two
fixed points coincide.
\end{prop}
\begin{proof}
Suppose first that the pair is strongly doubly reversible, and let $J$ be a common
involutory reverser. By Lemma~\ref{loxo1}, $J$ interchanges the two fixed points of
each element. If $p$ is a common fixed point and $p_1,p_2$ are the other fixed points
of $A,B$, respectively, then $J(p)=p_1=p_2$. Thus the fixed-point pairs coincide.

Conversely, suppose $A$ and $B$ have the same fixed points $p,q$. Conjugating
simultaneously in $\Sp(1,1)$, assume these are $o,\infty$ and write
\[
A_1=\begin{pmatrix}re^{i\theta}&0\\0&r^{-1}e^{i\theta}\end{pmatrix},\qquad
B_1=\begin{pmatrix}\mu&0\\0&\bar\mu^{-1}\end{pmatrix}.
\]
Choose $b\in\C^*$ so that, when $\mu=c_1+c_2j$, one has
$\operatorname{Re}(b\bar c_2)=0$, and set
\[
D=\begin{pmatrix}0&bj\\ \bar b^{-1}j&0\end{pmatrix}.
\]
Then $D^2=-I$ and a direct calculation gives
$DA_1D^{-1}=A_1^{-1}$ and $DB_1D^{-1}=B_1^{-1}$. Conjugating $D$ back gives a
common involutory reverser of $A$ and $B$ in $\PSp(1,1)$.
\end{proof}

\subsection{Pairs without a common fixed point}
Recall that for three distinct boundary points $p_1,p_2,p_3$, with lifts
$\mathbf p_1,\mathbf p_2,\mathbf p_3$, the quaternionic Cartan angular invariant is
\[
\mathbb A(p_1,p_2,p_3)=\arccos\left(
\frac{\operatorname{Re}(-H(p_1,p_2,p_3))}{|H(p_1,p_2,p_3)|}\right),
\]
where
\[
H(p_1,p_2,p_3)=\langle\mathbf p_1,\mathbf p_2\rangle
\langle\mathbf p_2,\mathbf p_3\rangle
\langle\mathbf p_3,\mathbf p_1\rangle.
\]
This invariant does not distinguish configurations in quaternionic hyperbolic
dimension one. Indeed, in the Siegel model a finite boundary point has the form
$(x,1)^t$ with $x\in\operatorname{Im}\H\cong\mathbb R^3$. For three finite
boundary points $x,y,z$, the real part of the Hermitian triple product is, up to a
nonzero real normalization, the scalar triple product of the three difference vectors
$x-y$, $y-z$, and $z-x$. Since these vectors sum to zero, that scalar triple product
vanishes. Thus
\[
\mathbb A(p_1,p_2,p_3)=\frac\pi2
\]
for every triple of distinct points of $\partial\h^1$; triples involving $\infty$
follow by invariance. Hence an equality or inequality involving only this angular
invariant cannot give a non-trivial criterion for strong double reversibility in
$\PSp(1,1)$.

We therefore turn directly to the matrix criterion, which records not only the exchange
of fixed points but also the rotational data required to reverse both elements.

\subsection{Quantitative description of strongly doubly reversible pairs in $\PSp(1,1)$}\label{qnt}

Throughout this subsection let
\[
A=\begin{pmatrix}re^{i\theta}&0\\0&r^{-1}e^{i\theta}\end{pmatrix}\in\Sp(1,1),
\qquad 0<r<1,\qquad
\theta\in(0,\pi)\setminus\{\pi/2\}.
\]
If $C\in\Sp(1,1)$ is a lift of an element of $\PSp(1,1)$ reversing $A$, then
$CAC^{-1}=\varepsilon A^{-1}$ with $\varepsilon\in\{\pm1\}$. Since
$\theta\neq\pi/2$, Lemma~\ref{invv}(iii) excludes $\varepsilon=-1$, and
Lemma~\ref{invv}(i) then shows that $C=\pm C_t$, where
\[
C_t=\begin{pmatrix}0&tj\\ \bar t^{-1}j&0\end{pmatrix},\qquad t\in\C^*.
\]
Conversely, $C_tAC_t^{-1}=A^{-1}$ and $C_t^2=-I$ for every $t\in\C^*$. Thus the
reversers of $A$ in $\PSp(1,1)$ are exactly the images of the $C_t$, and all of
them are involutions.

Let
\[
B=\begin{pmatrix}a&b\\c&d\end{pmatrix}\in\Sp(1,1),
\]
and write each entry as $h=h_1+h_2j$ with $h_1,h_2\in\C$, for $h=a,b,c,d$.
Since $B\in\Sp(1,1)$,
\[
B^{-1}=H_1B^*H_1=\begin{pmatrix}\bar d&\bar b\\ \bar c&\bar a\end{pmatrix}.
\]
For $\delta\in\{\pm1\}$ consider the system of equations
\begin{equation}\label{eq:qnt-core}
\begin{gathered}
(1-\delta)\,a_1=0,\qquad (1-\delta)\,d_1=0,\\
t\,\overline{a_2}+\delta\,\bar t\,a_2=0,\qquad
t\,\overline{d_2}+\delta\,\bar t\,d_2=0,\\
b_1=\delta\,|t|^2c_1,\qquad
b_2=-\delta\,t^2\,\overline{c_2}.
\end{gathered}
\end{equation}
For $\delta=1$ the first line is empty and the second line reads
$\operatorname{Re}(\bar t a_2)=\operatorname{Re}(\bar t d_2)=0$. For
$\delta=-1$ the system requires $a_1=d_1=0$ and
$\operatorname{Im}(\bar t a_2)=\operatorname{Im}(\bar t d_2)=0$.

\begin{lemma}\label{qnt-reduction}
Let $t\in\C^*$ and $\delta\in\{\pm1\}$. Then
$C_tBC_t^{-1}=\delta B^{-1}$ if and only if $(t,\delta)$ satisfies
\eqref{eq:qnt-core}.
\end{lemma}

\begin{proof}
The relation $C_tB=\delta B^{-1}C_t$ is equivalent to the four entrywise equations
\[
tjc=\delta\,\bar b\,\bar t^{-1}j,\qquad
tjd=\delta\,\bar d\,tj,\qquad
\bar t^{-1}ja=\delta\,\bar a\,\bar t^{-1}j,\qquad
\bar t^{-1}jb=\delta\,\bar c\,tj.
\]
Using $jz=\bar z j$ for $z\in\C$ and $j^2=-1$, we have
$jh=\overline{h_1}\,j-\overline{h_2}$ and $\bar h=\overline{h_1}-h_2j$. Comparing
the $\C$-components and the $\C j$-components of each equation gives the following.
\begin{itemize}
\item[(i)] From the $(1,2)$ entry: $t\overline{d_1}=\delta\,\overline{d_1}\,t$ and
$-t\overline{d_2}=\delta\,\bar t d_2$, that is, $(1-\delta)d_1=0$ and
$t\overline{d_2}+\delta\bar t d_2=0$.
\item[(ii)]  From the $(2,1)$ entry: $(1-\delta)a_1=0$ and, after multiplying by
$|t|^2$, $t\overline{a_2}+\delta\bar t a_2=0$.
\item[(iii)] From the $(1,1)$ entry: $\overline{b_1}=\delta|t|^2\overline{c_1}$ and
$b_2=-\delta t^2\overline{c_2}$.
\item[(iv)]  From the $(2,2)$ entry: the same two relations as from the $(1,1)$ entry.
\end{itemize}
Conversely, these relations imply the four entrywise equations.
\end{proof}

\begin{theorem}\label{qd}
With $A$ as above, $(A,B)$ is strongly doubly reversible in $\PSp(1,1)$ if and
only if there exist $t\in\C^*$ and $\delta\in\{\pm1\}$ satisfying
\eqref{eq:qnt-core}.
\end{theorem}

\begin{proof}
Suppose $(A,B)$ is strongly doubly reversible, and let $\iota\in\PSp(1,1)$ be a
common involutory reverser. Since $\iota$ reverses $A$, it has a lift $C_t$ for
some $t\in\C^*$. Since $\iota$ also reverses $B$, we have
$C_tBC_t^{-1}=\delta B^{-1}$ for some $\delta\in\{\pm1\}$. By
Lemma~\ref{qnt-reduction}, $(t,\delta)$ satisfies \eqref{eq:qnt-core}.

Conversely, if $(t,\delta)$ satisfies \eqref{eq:qnt-core}, then
Lemma~\ref{qnt-reduction} gives $C_tBC_t^{-1}=\delta B^{-1}$. We also have
$C_tAC_t^{-1}=A^{-1}$ and $C_t^2=-I$. Since $C_t$ is antidiagonal,
$C_t\neq\pm I$. Hence the image of $C_t$ in $\PSp(1,1)$ is an involution
reversing both $A$ and $B$.
\end{proof}

\begin{remark}
The restriction $\theta\in(0,\pi)\setminus\{\pi/2\}$ is exactly what makes every
reverser of $A$ of the form $C_t$. For $\theta\in\{0,\pi/2,\pi\}$,
Lemma~\ref{invv} produces further reversers, and Theorem~\ref{qd} does not apply.
For hyperbolic $B$, Theorem~\ref{thm:n1equiv} still applies in those cases.
\end{remark}

When the relevant complex components do not vanish, the criterion can be written
without the parameter $t$.

\begin{corollary}\label{qd-nondegenerate}
With $A$ as above, assume in addition that $a_2d_2c_1\neq0$. Then $(A,B)$ is
strongly doubly reversible in $\PSp(1,1)$ if and only if the following four
conditions hold:
\begin{enumerate}
\item[\textup{(i)}] $a_2\in\mathbb R\,d_2$;
\item[\textup{(ii)}] $b_1/c_1\in\mathbb R\setminus\{0\}$;
\item[\textup{(iii)}] $|a_2|^2\,b_2=\bigl|b_1/c_1\bigr|\,a_2^2\,\overline{c_2}$;
\item[\textup{(iv)}] if $b_1/c_1<0$, then $a_1=d_1=0$.
\end{enumerate}
In that case the sign in \eqref{eq:qnt-core} is
$\delta=\operatorname{sgn}(b_1/c_1)$.
\end{corollary}

\begin{proof}
Suppose $(t,\delta)$ satisfies \eqref{eq:qnt-core}.

\emph{Condition (ii) and the sign.} Since $c_1\neq0$, the relation
$b_1=\delta|t|^2c_1$ gives $b_1/c_1=\delta|t|^2$. So $b_1/c_1$ is a nonzero
real number, $\delta=\operatorname{sgn}(b_1/c_1)$, and
$|t|^2=|b_1/c_1|$.

\emph{Condition (i).} Since $a_2\neq0$, the equation
$t\overline{a_2}+\delta\bar t a_2=0$ says that $t\in i\,a_2\mathbb R^*$ if
$\delta=1$, and $t\in a_2\mathbb R^*$ if $\delta=-1$. The same holds with $d_2$
in place of $a_2$, since $d_2\neq0$. Both descriptions of $t$ are compatible
exactly when $a_2\in\mathbb R d_2$.

\emph{Condition (iii).} In either case write $t=s\,i\,a_2$ (if $\delta=1$) or
$t=s\,a_2$ (if $\delta=-1$), with $s\in\mathbb R^*$. Then $t^2=-\delta s^2a_2^2$
and $|t|^2=s^2|a_2|^2$, so $s^2=|b_1/c_1|/|a_2|^2$. The last equation of
\eqref{eq:qnt-core} becomes
\[
b_2=-\delta t^2\overline{c_2}=s^2a_2^2\,\overline{c_2}
=\frac{|b_1/c_1|}{|a_2|^2}\,a_2^2\,\overline{c_2}.
\]

\emph{Condition (iv).} If $b_1/c_1<0$, then $\delta=-1$, and the first line of
\eqref{eq:qnt-core} gives $a_1=d_1=0$.

Conversely, assume (i)--(iv). Put $\delta=\operatorname{sgn}(b_1/c_1)$ and
$s=|b_1/c_1|^{1/2}/|a_2|$. Let $t=s\,i\,a_2$ if $\delta=1$, and $t=s\,a_2$ if
$\delta=-1$. By (iv), the first line of \eqref{eq:qnt-core} holds. By (i), the
second line holds for both $a_2$ and $d_2$. By construction,
$\delta|t|^2c_1=|b_1/c_1|\,\delta\, c_1=b_1$. By (iii),
$-\delta t^2\overline{c_2}=s^2a_2^2\overline{c_2}=b_2$. Theorem~\ref{qd} now
completes the proof.
\end{proof}
\section*{Declarations} 
\subsection*{Authors’ Contributions}
       All authors contributed equally to this work.    
       \subsection*{Funding Declaration}
       Not applicable. 
	
	\subsection*{Ethical Approval} Not applicable.

\subsection*{Conflict of interest}There is no conflict of interest.

\subsection*{ Availability of data and materials} Not applicable.

\subsection*{AI Declaration}
The original manuscript was prepared without the use of any AI tools. In preparing the revised version, AI tools were used solely for proof verification and language editing.

\end{document}